# NONLINEAR SEQUENTIAL DESIGNS FOR LOGISTIC ITEM RESPONSE THEORY MODELS WITH APPLICATIONS TO COMPUTERIZED ADAPTIVE TESTS[1]

By Hua-Hua Chang and Zhiliang Ying

*University of Illinois at Urbana–Champaign and Columbia University*

Computerized adaptive testing is becoming increasingly popular due to advancement of modern computer technology. It differs from the conventional standardized testing in that the selection of test items is tailored to individual examinee's ability level. Arising from this selection strategy is a nonlinear sequential design problem. We study, in this paper, the sequential design problem in the context of the logistic item response theory models. We show that the adaptive design obtained by maximizing the item information leads to a consistent and asymptotically normal ability estimator in the case of the Rasch model. Modifications to the maximum information approach are proposed for the two- and three-parameter logistic models. Similar asymptotic properties are established for the modified designs and the resulting estimator. Examples are also given in the case of the two-parameter logistic model to show that without such modifications, the maximum likelihood estimator of the ability parameter may not be consistent.

**1. Introduction.** Computerized adaptive testing (CAT) is becoming increasingly popular due to advancement of modern computer technology. The concept of adaptive testing was originally conceived by Lord (1971) in his attempt to utilize the stochastic approximation algorithm of Robbins and Monro (1951) for designing more efficient tests. Major advances were carried out and documented in Owen (1975), Weiss (1976) and Wainer (2000). A distinctive feature of adaptive testing is to tailor test items (questions) to each examinee's ability level, so that able examinees can avoid doing too many

Received April 2008.

[1]Supported in part by the National Science Foundation, the National Security Agency and ETS research allocation PJ 79427.

easy items and less able examinees can avoid doing too many difficult items. Specifically, if the examinee answers a question correctly (incorrectly), then the next question administered to him/her will tend to be easier (more difficult). Through such an adaptive approach, questions with their difficulty levels suitable to a specific examinee are likely to be allocated. In consequence, examinees are challenged but not discouraged, leading to their ability levels being measured more accurately with the same or fewer number of items than using the conventional tests. Rapid development of computer technology has made adaptive testing a very promising option and, to a certain extent, the future of standardized tests. For example, computerized adaptive tests have already been implemented in GRE, the Graduate Record Examination, and GMAT, the Graduate Management Admission Test.

Both theoretical and implementational aspects of adaptive testing rely heavily on the item response theory (IRT) models, which relate examinees' ability levels to their responses to test items. Suppose that an examinee's ability level is characterized by a single parameter $\theta$. A basic assumption of the IRT is that for a given item the probability of producing a correct answer depends only on examinee's ability parameter $\theta$. The resulting probability curve, as $\theta$ varies, is known as the item characteristic curve (ICC) of the given item. Different parametrizations of the ICC curve lead to different IRT models.

Rasch (1960) proposed using the family of shifted logistic functions, $\exp(\theta - b)/(1 + \exp(\theta - b))$, to model the ICC. Here, $b$ determines the position of the ICC along the ability scale and is known as the item difficulty parameter. Exponent $\theta - b$ may be replaced by $1.7(\theta - b)$ to bring the curve closer to the standard normal distribution function. The latter will not be used in this paper, however, for mathematical simplicity. Let $Y$ denote an examinee's response, with values 1 indicating a correct answer and 0 an incorrect answer, to an item whose ICC follows the Rasch model with difficulty $b$. Then,

$$(1.1) \qquad P(Y = 1 | \theta) = \frac{e^{\theta - b}}{1 + e^{\theta - b}},$$

where $\theta$ denotes the ability level of the examinee.

A more general model, which includes the Rasch model as a special case, is the so-called three-parameter logistic (3-PL) model, whose ICC is defined by

$$(1.2) \qquad P(Y = 1 | \theta) = c + (1 - c)\frac{e^{a(\theta - b)}}{1 + e^{a(\theta - b)}},$$

where $Y$, $\theta$ and $b$ have the same interpretations as those in (1.1) and where the additional item parameters $c$ and $a$ measure, respectively, the degree of guessing and the discriminating power [see Lord (1980), Hambleton and Swaminathan (1985)]. The Rasch model corresponds to the situation in which there is no guessing ($c \equiv 0$) and all items have the same discriminating power ($a \equiv 1$ when properly scaled). An intermediate model is

$$(1.3) \qquad P(Y = 1 | \theta) = \frac{e^{a(\theta - b)}}{1 + e^{a(\theta - b)}},$$

which is known as the two-parameter logistic (2-PL) model.

The conventional IRT model-based design of a test is the advance selection of a set of $n$ items whose parameters have been precalibrated (known). For each examinee, there are $n$ responses, say $Y_1, \ldots, Y_n$, to the $n$ items. Point and interval estimation of $\theta$ for the examinee can then be obtained by, for example, maximizing the likelihood function of $\theta$ with $Y_1, \ldots, Y_n$ and calculating the observed Fisher information, or by other methods that can be found in statistical literature. Lord (1980) contains detailed descriptions of relevant statistical inference procedures and theory thereof.

The main focus of the present investigation is on the IRT model-based adaptive design of computerized tests. An adaptive test differs from a conventional test in that the assignment of the test items are performed sequentially, with selection of each item depending on the responses of the examinee to the preceding items. More specifically, let $A$ be the item bank from which items may be selected and assigned to the examinee. Suppose that $k-1$ items, $\alpha_1, \ldots, \alpha_{k-1} \in A$ have already been selected and that the responses from the examinees are $Y_1, \ldots, Y_{k-1}$. The selection of the $k$th item, $\alpha_k$, will be based on the previous items, $\alpha_1, \ldots, \alpha_{k-1}$ as well as the responses $Y_1, \ldots, Y_{k-1}$. Arising from this formulation are three aspects that may be studied: (1) Design of an adaptive rule for selection of test items $\alpha_1, \alpha_2, \ldots$, (2) sequential estimation of ability parameter $\theta$ at each stage, and (3) properties of the adaptive design and the resulting estimator.

Lord (1980), Chapter 10, argued that, for a given examinee, the items should be selected to maximize the Fisher information. Let $P_\alpha(\theta)$ be the probability that an examinee with ability $\theta$ answers item $\alpha$ correctly. The Fisher information function (of $\theta$) for $\alpha$ is simply

$$(1.4) \qquad I_\alpha(\theta) = \left[ \frac{\partial P_\alpha}{\partial \theta}(\theta) \right]^2 \Big/ P_\alpha(\theta) Q_\alpha(\theta),$$

where $Q_\alpha(\theta) = 1 - P_\alpha(\theta)$. If $\theta$ were known, then the optimal choice, according to Lord, is the one that maximizes $I_\alpha(\theta)$. Although in reality we do not know $\theta$, the sequential approach allows us to use the current estimate of $\theta$ in deciding the next choice of $\alpha$. Our results, to be presented in this paper, indicate that, for the Rasch model, such an approach leads to an asymptotically optimal design and that, for the two-parameter and three-parameter logistic model, the approach does not in general lead to an optimal design. In fact, the procedure needs to be modified in order to produce a reasonable design. Note that, throughout the paper, the term *optimal* is referred to that the adaptive design leads to a consistent and asymptotically normal ability estimator.

Despite the increased prominence of CAT in standardized testing, in-depth statistical analysis has yet to be developed. The present paper is aimed at providing some basic results in certain idealized situations. It is organized as follows. The Rasch model is studied in Section 2 in the context of the adaptive design and maximum likelihood estimation. It is shown there that the maximum Fisher information-based sequential design, in conjunction with updating maximum likelihood recursion, is asymptotically optimal, and the resulting maximum likelihood estimator is consistent and asymptotically normal. In Section 3, a modification to the maximum Fisher information-based design for the two-parameter logistic model is proposed, and the resulting maximum likelihood estimator is shown to be consistent and asymptotically normal. A counterexample is also given to illustrate the necessity of the proposed modification. Treatment of the general three-parameter logistic model is given in Section 4, where, in addition to modifying the maximum Fisher information design, we also propose an approximation to the maximum likelihood estimating equation. The usual large sample properties are established accordingly. Discussions and some concluding remarks are given in Section 5.

## 2. Information-based adaptive design for the Rasch model.

Recall that, under the Rasch model, the probability of answering an item correctly by an examinee with ability parameter $\theta$ is $\exp(\theta - b)/[1 + \exp(\theta - b)]$, where $b$ is the item parameter representing the difficulty level. From (1.4), the Fisher information of the item can be written as

$$(2.1) \qquad I(\theta|b) = \frac{e^{\theta-b}}{(1+e^{\theta-b})^2}.$$

For a given examinee, $I_b(\theta)$ attains its maximum value $1/4$ at $b = \theta$. Therefore, the optimal design is to select items with difficulty parameter $b = \theta$. Since $\theta$ is unknown, successive approximations to the optimal design will be needed.

A general recursive algorithm known as the stochastic approximation for approximating optimal design points was first proposed by Robbins and Monro (1951). Lord (1971) discovered use of the stochastic approximation in developing adaptive (tailored) tests. Wu (1985) introduced a maximum likelihood recursion as an alternative to the stochastic approximation when the underlying response curve is of the logistic form. He further showed, through extensive simulation studies, that his maximum likelihood recursion improves efficiency over the stochastic approximation when the sample size is moderate.

In this section, we first consider an idealized setting for CAT in which available items at each stage exhaust all difficulty levels. In other words, for every $b$, an item with ICC $\{\exp(\theta - b)/[1 + \exp(\theta - b)], \theta \in R\}$ can be administered to the examinee. We will then consider more realistic situations for which available items are limited, so that we can at best choose items that are closest to the idealized optimal ones. Results for the idealized CAT will be developed and then extended for the more realistic situations.

For the idealized CAT, the sequential design based on maximizing the Fisher information and updating maximum likelihood estimators consists of the following steps:

1. *Initialization.* Specify the difficulty level, say $b_1$, of the initial item. If the examinee's response is correct (i.e., $Y_1 = 1$), then choose the succeeding items with increasing difficulty parameters $(b_1 \leq) b_2 \leq b_3 \cdots \leq b_{k_0}$, where $k_0 = \inf\{j : Y_j = 0\}$ is the first time an incorrect response occurs. On the other hand, if the response to the first item is incorrect, then select the succeeding items with decreasing difficulty parameters $(b_1 \geq) b_2 \geq b_3 \cdots \geq b_{k_0}$, where $k_0 = \inf\{j : Y_j = 1\}$.

2. *Estimation.* For each $k \geq k_0$, define $\hat{\theta}_k$ by solving the maximum likelihood estimating equation
$$\sum_{i=1}^{k} \left( Y_i - \frac{e^{\theta - b_i}}{1 + e^{\theta - b_i}} \right) = 0.$$
Since the response sequence $\{Y_1, \ldots, Y_k\}$ contains both 0 and 1, $\hat{\theta}_k$ is uniquely and well defined.

3. *Design.* After $k(\geq k_0)$ items are administered and $\hat{\theta}_k$ is obtained, select the next item by setting $b_{k+1} = \hat{\theta}_k$. Note that this selection is simply the idealized optimal design, but with unknown parameter $\theta$ being replaced by its most recent estimator.

The preceding adaptive testing procedure was proposed and discussed in Lord (1971, 1980). It was also studied in the context of sequential optimal design in Wu (1985), where its connection to Robbins and Monro's stochastic approximation algorithm was found. Ying and Wu (1997) established an asymptotic theory for a class of sequential design problems. The next theorem shows that the sequential estimator $\hat{\theta}_n$ is consistent and asymptotically normal. It entails that the adaptive design is asymptotically efficient.

THEOREM 1. *Let $\{\hat{\theta}_k\}$ be the sequential estimators specified by steps 1–3 for the Rasch model. Then, as $n \to \infty$, $\hat{\theta}_n \to \theta$ a.s. and $\sqrt{n/4}(\hat{\theta}_n - \theta) \to N(0, 1)$. Furthermore, $4I_n(\hat{\theta})/n \to 1$ a.s., where $I_n(\theta) = \sum_{i=1}^{n} \exp(\theta - b_i)/(1 + \exp(\theta - b_i))^2$ is the observed Fisher information.*

The asymptotic variance for $\hat{\theta}_n$ is $4/n$, which is exactly the inverse of the Fisher information if all the $n$ items are chosen optimally (i.e., $b_i \equiv \theta$). Thus, the estimator $\hat{\theta}_n$ is asymptotically optimal. However, under the more realistic situation in which the item bank has limited capacity, that is $b_k$ can only be chosen from a set of discrete values, then the consistency and asymptotic normality for $\hat{\theta}_k$ still hold, but the asymptotic variance needs to be replaced by the inverse of the Fisher information.

Theorem 1 is implied by the more general result given by Theorem 2. It can also be inferred from Ying and Wu (1997), Theorem 1. Proof of Theorem 2 uses the so-called local convergence theorem for martingale sequences.

**3. The two-parameter logistic model.** Recall that the two-parameter logistic model is an extension of the Rasch model to include a second item parameter $a$, which represents the discriminating power of the item. Under this model, an examinee with ability $\theta$ answers an item, specified by $a$ and $b$, correctly with probability $e^{a(\theta-b)}/[1 + e^{a(\theta-b)}]$ in (1.3).

The Fisher information function for an item specified by $a$ and $b$ may be expressed as

$$(3.1) \qquad I(\theta|a,b) = a^2 \frac{e^{a(\theta-b)}}{[1 + e^{a(\theta-b)}]^2}.$$

If $a$ and $b$ are unrestricted, then the information-based optimal design problem is singular because

$$(3.2) \qquad \max_{a,b} I(\theta|a,b) \left( = \max_a \max_b I(\theta|a,b) \right) = \max_a \frac{a^2}{4} = \infty.$$

From (3.2), the optimal design appears to be $b = \theta$ and $a = \infty$. But this will be extremely unstable since, for any $b \neq \theta$,

$$\lim_{a \to \infty} I(\theta|a,b) = 0.$$

One way to avoid such singularity is to restrict the item pool so that parameter $a$ will fall into a compact interval in $(0, \infty)$.

Analogous to the adaptive design for the Rasch model, we introduce a similar design for the two-parameter logistic model. However, to avoid the singularity, we shall put a restriction on the discrimination parameter $a$. Specifically, let $0 < m < M < \infty$ be fixed in advance, and assume $a \in [m, M]$.

1. *Initialization.* Select the initial coin (item) with parameters $a_1$ and $b_1$. Reasonable choice for them can be made from the prior information about the population. If the outcome of the first toss $Y_1$ is 1 (head), then choose the next $k_0$ coins with increasing difficulty parameters $(b_1 \leq) b_2 \leq \cdots \leq b_{k_0}$, where $k_0 = \inf\{j : Y_j = 0\}$ is again the first time a tail occurs. If the first toss is a tail, then choose $(b_1 \geq)b_1 \geq \cdots \geq b_{k_0}$ with $k_0$ being the first head. The $a$-parameters must satisfy $m \leq a_j \leq M, j = 1, \ldots, k_0$ but can be arbitrary otherwise.

2. *Estimation.* For each $k \geq k_0$, define $\hat{\theta}_k$, the maximum likelihood estimator, as the unique solution to

$$(3.3) \qquad \sum_{i=1}^{k} a_i \left( Y_i - \frac{e^{a_i(\theta - b_i)}}{1 + e^{a_i(\theta - b_i)}} \right) = 0.$$

Note that the left-hand side of (3.3) is a strictly decreasing function with values ranging from $\sum_{i=1}^{k} a_i(1 - Y_i) < 0$ to $\sum_{i=1}^{k} a_i Y_i > 0$.

3. *Design.* After $\hat{\theta}_k$ is defined, set $b_{k+1} = \hat{\theta}_k$ and $a_{k+1}$ to be a number in $[m, M]$. The choice for $a_{k+1}$ can depend on data collected up to the current stage. The next selection will be the coin (item) with parameters $a_{k+1}$ and $b_{k+1}$.

The preceding sequential design is not optimal, not even asymptotically. It is based on a suboptimal design that maximizes the Fisher information over $b$ with $a$ being fixed. Such an approach is intuitively sensible, because the adaptive test is to match the difficulty level of test items with examinee's ability and parameter $b$ represents the item difficulty. Obviously, it does not touch upon selection of the discrimination parameter $a$, which involves more complex issues [see Chang and Ying (1999) and Chang, Qian and Ying (2001)].

THEOREM 2. *Under the preceding sequential design for the two-parameter logistic model, $\hat{\theta}_n \to \theta$ a.s. as $n \to \infty$. In addition, suppose the choice of $a_j$ satisfies $\sum_{i=1}^{n} a_i^2 / v_n \to_p 1$, as $n \to \infty$, for some nonrandom sequence $v_n$. Then,*

$$(3.4) \qquad \sqrt{\sum_{i=1}^{n} a_i^2}(\hat{\theta}_n - \theta) \to_{\mathcal{L}} N(0, 1).$$

*The normalizing factor $\sqrt{\sum_{i=1}^{n} a_i^2}$ in (3.4) may be replaced by $\sqrt{I^{(n)}(\hat{\theta}_n)}$ or $\sqrt{I^{(n)}(\theta)}$, where*

$$(3.5) \qquad I^{(n)}(\theta) = \sum_{i=1}^{n} a_i^2 \frac{e^{a_i(\theta - b_i)}}{[1 + e^{a_i(\theta - b_i)}]^2}$$

*is the observed Fisher information.*

REMARK 1. As we stated earlier, the solution to the optimal design problem of maximizing the Fisher information is singular, in that the discrimination parameter will reach $\infty$. The remedial measure taken here is to restrict this parameter to a compact interval. Next we construct an example to show that if the $a_j$ are not bounded, it is possible that the resulting estimator $\hat{\theta}_n$ may not even be consistent.

EXAMPLE 1. Suppose we follow the same sequential design as described at the beginning of this section, but with $a_k = k^3$ instead of confining the $a_k$ to a compact interval. Suppose, in addition, that the initial value is taken to be $\hat{\theta}_0 < \theta - 1 - \pi^2/6$. If $Y_1 = \cdots = Y_j = 0$, then the subsequent $\hat{\theta}_k, 1 \le k \le j$, will be chosen in decreasing order, so that $\hat{\theta}_1 = \hat{\theta}_0 - \varepsilon_0, \ldots, \hat{\theta}_j = \hat{\theta}_{j-1} - \varepsilon_0$, where $\varepsilon_0 > 0$ is a prespecified constant. Let $n_0$ be a large integer, so that the following conditions are satisfied:

$$\text{(3.6)} \qquad \sum_{k=n_0+1}^{\infty} \frac{k^3}{1+e^k} < \frac{1}{3},$$

$$\text{(3.7)} \qquad \frac{(n_0+1)^3}{1+e^{(n_0+1)^3 \varepsilon_0}} < \frac{1}{6},$$

$$\text{(3.8)} \qquad 3(n_0+1)^3 < \sum_{k=1}^{n_0} k^3.$$

Define event $A = \{Y_k = 0, k \le n_0 \text{ and } Y_k = 1, k \ge n_0 + 1\}$. We prove below that $P(A) > 0$ and $\lim_{n \to \infty} \hat{\theta}_n < \theta - 1$ on $A$. Therefore, with such a design, $\hat{\theta}_n$ cannot be a consistent estimator of $\theta$. Intuitively, this can occur because movement of successive $\hat{\theta}_j$ is tied to the $a$-parameter. A large value of $a$ corresponds to a small movement size. The constructed example makes the $a$-parameters so large that the $\hat{\theta}_j$ can never move back, even if all the steps after $n_0$ are in right direction. Figure 1 shows graphically two sequences of $\hat{\theta}_j$, one converges to the $\theta$ and the other does not.

REMARK 2. The constraint that the discrimination parameters $a_j$ are bounded away from 0 is also needed. To see this, suppose we set $a_j = \frac{1}{j}$. Then, the total Fisher information for a test of length $n$ is bounded by

$$\frac{1}{4} \sum_{j=1}^{n} a_j^2 < \frac{1}{4} \sum_{j=1}^{\infty} \frac{1}{j^2} < \infty.$$

In view of this, it is straightforward that the resulting maximum likelihood estimator $\hat{\theta}_n$ will not converge to $\theta$.

REMARK 3. It was pointed out by the Associate Editor that an item with a large value of $a$-parameter could be very uninformative if knowledge about $\theta$ is poor, and also that a natural way to increase efficiency is to use items with small $a$-parameter values in early stages and use items with large $a$-parameter values in later stages. Indeed, such an approach could lead to, among other things, substantial efficiency improvement. Figure 2 gives efficiency comparison in terms of mean squared errors using ascending and descending $a$-parameter values. For more details and other related issues in practical settings, we refer to Chang and Ying (1996, 1999). It is worth noting that, if items with $a = \infty$ were available, then one could design a scheme that approaches the true $\theta$ exponentially fast, though such a scheme is likely to be different from the maximum likelihood estimation.

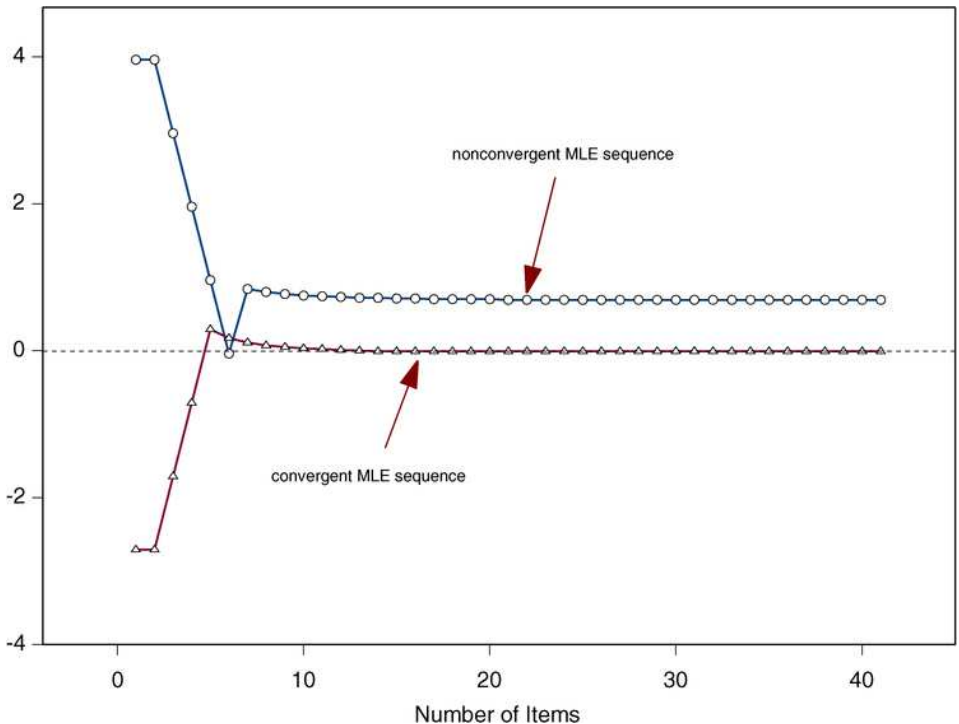

Fig. 1. *Examples of convergence and nonconvergence.*

REMARK 4. As pointed by one of the reviewers, setting $a \leq M$ is reasonable, because no item-writer has ever been able to write a sequence of items with $a$-parameters tending to infinity. Also, any reasonable item bank would only include items with $a$-parameters bounded away from 0. If we assume the item bank contains all pairs $(a, b)$ in $[m, M]$ by $(-\infty, \infty)$, Chang and Ying (1999) proposed the $a$-stratified method with an objective to limit the exposure on any given item by using that item at the most advantageous point in testing. The $a$-stratified method attempts to use less discriminating items early in the test, when estimation is least precise, and save highly discriminating items until later stages, when finer gradations of estimation are required. One of the advantages of using the $a$-stratified method is that it attempts to equalize the item exposure rates for all the items in the pool.

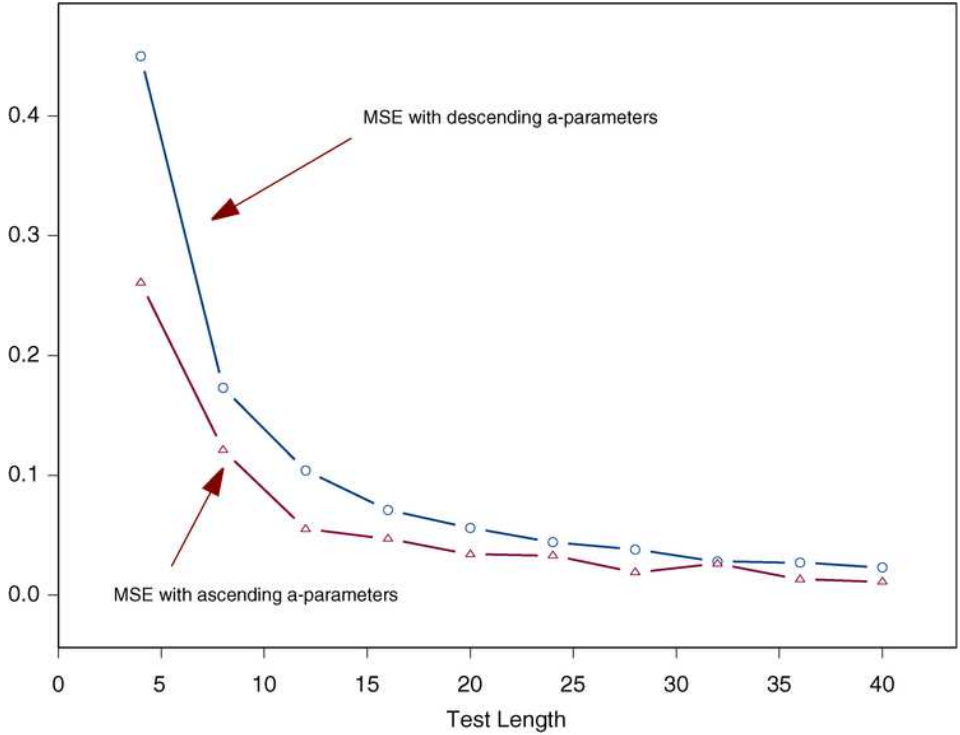

Fig. 2. *Mean squared errors under ascending and descending a-parameter designs.*

PROOF OF THEOREM 2. The main line of the proof consists of the following four steps. First, we show that the observed Fisher information goes to infinity as $n \to \infty$. The second step is to show that the design leads to bounded maximum likelihood estimators $\hat{\theta}_n$. From the boundedness follows the consistency. The last step is to show the asymptotic normality.

Throughout the proof, we shall let $G(t) = e^t/(1+e^t)$ and $\bar{G}(t) = 1 - G(t)$. Define $\sigma$-filtration $\mathcal{F}_k = \sigma\{Y_j, \hat{\theta}_j, a_{j+1}, j \leq k\}$, $k \geq 0$. Then, conditioning on $\mathcal{F}_{k-1}, Y_k$ is a Bernoulli random variable with success probability $G(a_k(\theta - b_k))$. Thus, $\{Y_k - G(a_k(\theta - b_k))\}$ is a martingale difference sequence with respect to $\{\mathcal{F}_k\}$. Since $a_k \in \mathcal{F}_{k-1}$ is predictable, $a_k[Y_k - G(a_k(\theta - b_k))]$ are again martingale differences with conditional variances

$$\text{Var}\{a_k[Y_k - G(a_k(\theta - b_k))] | \mathcal{F}_{k-1}\} = a_k^2 G(a_k(\theta - b_k))\bar{G}(a_k(\theta - b_k)).$$

Applying the martingale local convergence theorem of Chow (1965), Corollary 5, we have that

$$(3.9) \qquad \sum_{k=1}^{\infty} \frac{a_k[Y_k - G(a_k(\theta - b_k))]}{\sum_{j=1}^{k} a_j^2 G(a_j(\theta - b_j))\bar{G}(a_j(\theta - b_j))} \qquad \text{converges a.s.}$$

We first prove that

$$(3.10) \qquad P\left(\sum_{k=1}^{\infty} a_k^2 G(a_k(\theta - b_k))\bar{G}(a_k(\theta - b_k)) < \infty\right) = 0.$$

Let $A_1$ be the event that $\sum_{k=1}^{\infty} a_k^2 G(a_k(\theta - b_k))[1 - G(a_k(\theta - b_k))] < \infty$. Clearly, on $A_1$, $\lim_{n\to\infty} |b_n| = \infty$ or, equivalently, $\lim_{n\to\infty} |\hat{\theta}_n| = \infty$, recalling that $b_{n+1} = \hat{\theta}_n$ as the design requires it. From (3.9) and the monotonicity of the denominator sequence in (3.9), we have that $\sum_{k=1}^{\infty} a_k[Y_k - G(a_k(\theta - b_k))]$ converges on $A_1$. But $\sum_{k=1}^{n} a_k[Y_k - G(a_k(\hat{\theta}_n - b_k))] = 0$ for all $n$. So, on $A_1$,

$$\infty > \left|\sum_{k=1}^{\infty} a_k[Y_k - G(a_k(\theta - b_k))]\right|$$

$$(3.11) \qquad = \lim_{n\to\infty}\left|\sum_{k=1}^{n} a_k[Y_k - G(a_k(\theta - b_k))] - \sum_{k=1}^{n} a_k[Y_k - G(a_k(\hat{\theta}_n - b_k))]\right|$$

$$= \lim_{n\to\infty}\left|\sum_{k=1}^{n} a_k[G(a_k(\hat{\theta}_n - b_k)) - G(a_k(\theta - b_k))]\right|.$$

From (3.11) we claim that $\limsup \hat{\theta}_n < \infty$ on $A_1$. We prove this claim by contradiction. Suppose it is not true. Then there exists a subsequence $n_j$ such that $\hat{\theta}_{n_j} \to \infty$ and $\hat{\theta}_{n_j} \geq \hat{\theta}_k$ for all $k \leq n_j$. This implies the following:

$$(3.12) \qquad a_k[G(a_k(\hat{\theta}_{n_j} - b_k)) - G(a_k(\theta - b_k))] \geq 0 \qquad \text{for all } k \leq n_j;$$

$$(3.13) \qquad \text{for any constant } K, \qquad \#\{k \leq n_j : a_k(\theta - b_k) \leq -K\} \to \infty$$
$$\text{as } n_j \to \infty.$$

Combining (3.12) with (3.13), we know that (3.11) cannot be true. Thus, $\limsup \hat{\theta}_n < \infty$ on $A_1$. Likewise, we can show that $\liminf \hat{\theta}_n > -\infty$ on $A_1$. These two contradict a previous conclusion that $\limsup_{n\to\infty} |\hat{\theta}_n| = \infty$ on $A_1$ unless $A_1$ is a null event. Hence (3.10) holds.

From (3.9), (3.10) and the Kronecker lemma [Chow and Teicher (1988), page 114], it follows that

$$(3.14) \qquad \frac{\sum_{k=1}^{n} a_k[Y_k - G(a_k(\theta - b_k))]}{\sum_{k=1}^{n} a_k^2 G(a_k(\theta - b_k))\bar{G}(a_k(\theta - b_k))} \to 0 \qquad \text{a.s.}$$

Substituting the likelihood equation into (3.14), we get

$$(3.15) \qquad \frac{\sum_{k=1}^{n} a_k [G(a_k(\hat{\theta}_n - b_k)) - G(a_k(\theta - b_k))]}{\sum_{k=1}^{n} a_k^2 G(a_k(\theta - b_k)) \bar{G}(a_k(\theta - b_k))} \to 0 \qquad \text{a.s.,}$$

which certainly implies that

$$(3.16) \qquad \frac{1}{n} \sum_{k=1}^{n} a_k [G(a_k(\hat{\theta}_n - b_k)) - G(a_k(\theta - b_k))] \to 0 \qquad \text{a.s.}$$

Next, we show that $\limsup |\hat{\theta}_n| < \infty$ a.s. Suppose that this is not true and that, without loss of generality, there exists a subsequence $\{n_j\}$ such that $\hat{\theta}_{n_j} \uparrow \infty$ and $\hat{\theta}_{n_j} \geq \hat{\theta}_{k-1} = b_k$ for all $k \leq n_j$. Let $\delta_0 > 0$ be a fixed constant. Since $m \leq a_k \leq M$, we have

$$(3.17) \qquad a_k(\theta - b_k) \leq -\delta_0 \qquad \text{if } b_k \geq \theta + \frac{\delta_0}{m}.$$

But $G(a_k(\hat{\theta}_{n_j} - b_k)) \geq G(0) = \frac{1}{2}$ for all $k \leq n_j$, which, together with (3.17), implies that

$$(3.18) \qquad G(a_k(\hat{\theta}_{n_j} - b_k)) - G(a_k(\theta - b_k)) \geq \frac{1}{2} - G(-\delta_0) > 0$$
$$\text{if } b_k \geq \theta + \frac{\delta_0}{m}.$$

Thus, $\limsup_{j \to \infty} \#\{k \leq n_j : b_k \geq \theta + \frac{\delta_0}{m}\}/n_j = 0$, since we can otherwise select a subsequence of $n_j$ such that (3.16) does not hold. On the other hand,

$$(3.19) \qquad \begin{aligned} G(a_k(\hat{\theta}_n - b_k)) &- G(a_k(\theta - b_k)) \\ &= \frac{e^{a_k(\hat{\theta}_n - b_k)} - e^{a_k(\theta - b_k)}}{[1 + e^{a_k(\hat{\theta}_n - b_k)}][1 + e^{a_k(\theta - b_k)}]} \\ &= \frac{e^{a_k(\theta - b_k)}(e^{a_k(\hat{\theta}_n - \theta)} - 1)}{[1 + e^{a_k(\hat{\theta}_n - b_k)}][1 + e^{a_k(\theta - b_k)}]} \\ &= \frac{e^{a_k(\theta - b_k)}(1 + e^{-a_k(\hat{\theta}_n - \theta)})}{[e^{-a_k(\hat{\theta}_n - \theta)} + e^{a_k(\theta - b_k)}][1 + e^{a_k(\theta - b_k)}]}. \end{aligned}$$

Since $n_j \to \infty$, we have, in view of (3.19),

$$G(a_k(\hat{\theta}_{n_j} - b_k)) - G(a_k(\theta - b_k)) = (1 + o(1)) \frac{e^{a_k(\theta - b_k)}}{[o(1) + e^{a_k(\theta - b_k)}][1 + e^{a_k(\theta - b_k)}]},$$

which has the same order as

$$G(a_k(\theta - b_k)) \bar{G}(a_k(\theta - b_k)) = \frac{e^{a_k(\theta - b_k)}}{[1 + e^{a_k(\theta - b_k)}]^2}$$

for all $b_k \le \theta + \delta_0/m$. But, we know that $\#\{k \le n_j : b_k \le \theta + \delta_0/m\}/n_j \to 1$. So, (3.15) cannot hold along $n = n_j$. This contradiction proves that $\limsup |\hat{\theta}_n| < \infty$ a.s.

Now, by the mean value theorem, there exists $\hat{\theta}_n^*$ between $\theta$ and $\hat{\theta}_n$ such that

$$\sum_{k=1}^{n} a_k [G(a_k(\hat{\theta}_n - b_k)) - G(a_k(\theta - b_k))]$$

$$= \sum_{k=1}^{n} a_k^2 G(a_k(\hat{\theta}_n^* - b_k)) \bar{G}(a_k(\hat{\theta}_n^* - b_k))(\hat{\theta}_n - \theta).$$

Furthermore, $\liminf n^{-1} \sum_{k=1}^{n} a_k^2 G(a_k(\hat{\theta}_n^* - b_k)) \bar{G}(a_k(\hat{\theta}_n^* - b_k)) > 0$ since $\limsup |\hat{\theta}_n| < \infty$. Hence, (3.15) implies that $\hat{\theta}_n \to \theta$ a.s.

To prove the asymptotic normality, we follow the standard approach by taking the Taylor expansion; that is,

$$\begin{aligned}
0 &= \sum_{k=1}^{n} a_k [Y_k - G(a_k(\hat{\theta}_n - b_k))] \\
(3.20) \quad &= \sum_{k=1}^{n} a_k [Y_k - G(a_k(\theta - b_k))] \\
&\quad - \sum_{k=1}^{n} a_k^2 G(a_k(\hat{\theta}_n^* - b_k)) \bar{G}(a_k(\hat{\theta}_n^* - b_k))(\hat{\theta}_n - \theta),
\end{aligned}$$

where $\hat{\theta}_n^*$ is between $\hat{\theta}_n$ and $\theta$ and therefore converges to $\theta$ a.s. From (3.20), we have

$$\hat{\theta}_n - \theta = \left[ \sum_{k=1}^{n} a_k^2 G(a_k(\hat{\theta}_n^* - b_k)) \bar{G}(a_k(\hat{\theta}_n^* - b_k)) \right]^{-1} \sum_{k=1}^{n} a_k [Y_k - G(a_k(\theta - b_k))].$$

Since $\hat{\theta}_n^* \to \theta$ a.s. and $b_n \to \theta$ a.s., it follows that (3.4) holds if we can show that

$$(3.21) \quad \left( \sum_{k=1}^{n} a_k^2 \right)^{-1/2} \sum_{k=1}^{n} a_k [Y_k - G(a_k(\theta - b_k))] \to_{\mathcal{L}} N(0,1).$$

By the assumption, there is a nonrandom sequence $v_n \to \infty$ such that $\sum_{k=1}^{n} a_k^2 / v_n \to_p 1$. Thus, we can apply the martingale central limit theorem, as stated in Pollard [(1984), page 171] to get (3.21). Because $\hat{\theta} \to \theta$ a.s., we can easily see that $\sum_{k=1}^{n} a_k^2$ is asymptotically equivalent to $I^{(n)}(\hat{\theta}_n)$ as well as $I^{(n)}(\theta)$. $\square$

**4. The three-parameter logistic model and a modification to the maximum likelihood recursions.** The three-parameter logistic model, as specified by (1.2), extends the two-parameter model by including an additional parameter known as the guessing parameter. Recall that the ICC, in this case, is $c + (1-c)e^{a(\theta-b)}/[1 + e^{a(\theta-b)}]$. It is not difficult to see that, when $c > 0$, the family of probability distributions indexed by $\theta$ no longer forms an exponential family. Therefore, we expect that there will be extra technical difficulties to deal with.

For an item with parameters $a, b$ and $c$, the associated Fisher information function may be calculated using (1.4) to be

$$(4.1) \qquad I(\theta|a,b,c) = \frac{(1-c)a^2 e^{2a(\theta-b)}}{[c + e^{a(\theta-b)}][1 + e^{a(\theta-b)}]^2}.$$

For fixed $a$ and $c$, the Fisher information reaches its maximum when

$$(4.2) \qquad b = \theta - \frac{1}{a}\log\frac{1 + \sqrt{1+8c}}{2}$$

[see Lord (1980), page 152]. As the two examples indicated in Section 3, the discrimination parameter cannot be chosen arbitrarily because otherwise it may lead to inconsistency. It is also reasonable to put restrictions on selecting $c$, the guessing parameter. This is because, in view of (4.1), the Fisher information reaches the maximum if and only if $c = 0$. So if no constraint is put, then only those items with no guessing will be used. The design problem we shall be considering will only involve choice of $b$, with $a$ and $c$ being confined to certain reasonable regions.

In view of (4.2), we can select the optimal $b$ if $\theta$ is specified. The adaptive optimal design is then to replace $\theta$ by its current estimator. As we shall see, it turns out that the maximum likelihood estimating equation may have multiple roots. To avoid such a situation, we shall propose a modification, which is asymptotically equivalent to the likelihood estimating equation and has a unique root.

Suppose that the examinee has answered $n$ items, which are specified by $(a_k, b_k, c_k), k = 1, \ldots, n$, and the results are $Y_1, \ldots, Y_n$. Then, the maximum likelihood estimating equation for $\theta$ may be written as

$$(4.3) \qquad \sum_{k=1}^{n} \frac{a_k e^{a_k(\theta-b_k)}}{c_k + e^{a_k(\theta-b_k)}}\left[Y_k - c_k - (1-c_k)\frac{e^{a_k(\theta-b_k)}}{1 + e^{a_k(\theta-b_k)}}\right] = 0.$$

Unlike in the two-parameter logistic model, the left-hand side of (4.3) is not a monotone function of $\theta$. In fact, (4.3) may have multiple roots [Samejima (1973)]. On the other hand, when the choice of the difficulty parameter satisfies (4.2) ($\theta$ will be replaced by the current estimator), it is easy to see

that weights in (4.3)

$$\frac{a_k e^{a_k(\theta - b_k)}}{c_k + e^{a_k(\theta - b_k)}} \approx \frac{a_k(1 + \sqrt{1 + 8c_k})}{2c_k + 1 + \sqrt{1 + 8c_k}}.$$

Therefore, an approximation to (4.3) is

$$(4.4) \quad \sum_{k=1}^{n} \frac{a_k(1 + \sqrt{1 + 8c_k})}{2c_k + 1 + \sqrt{1 + 8c_k}} \left[ Y_k - c_k - (1 - c_k)\frac{e^{a_k(\theta - b_k)}}{1 + e^{a_k(\theta - b_k)}} \right] = 0,$$

which will be called approximate maximum likelihood estimating equation. It is obvious that the left-hand side of (4.4) is monotone decreasing in $\theta$. Therefore, the solution to (4.4), if it exists, will be unique. Notice also that the weights in (4.4) do not depend on the $b_k$.

An extension of the adaptive design procedure proposed in the preceding section to the three-parameter model is described below:

1. *Initialization.* In the same way as that for the two-parameter logistic model, choose the initial $k_0$ items so that $\{Y_i, i \leq k_0\}$ contains both 0 and 1.

2. *Selection of $\hat{\theta}_k$.* For each $k \geq k_0$, if

$$(4.5) \quad \sum_{i=1}^{k} \frac{a_i(1 + \sqrt{1 + 8c_i})}{2c_i + 1 + \sqrt{1 + 8c_i}} Y_i > \sum_{i=1}^{k} \frac{a_i(1 + \sqrt{1 + 8c_i})}{2c_i + 1 + \sqrt{1 + 8c_i}} c_i,$$

   then define $\hat{\theta}_k$ as the unique solution to (4.4). Otherwise, set $\hat{\theta}_k = r_k$, where $r_k \downarrow -\infty$ is a predetermined sequence.

3. *Design.* After selecting $\hat{\theta}_k$, set $b_{k+1} = \hat{\theta}_k$. Also, set $a_{k+1}$ and $c_{k+1}$ such that $a_{k+1} \in [m, M]$, and $c_{k+1} \leq 1 - \delta_0$, where $\delta_0 > 0$ is some constant.

REMARK 5. If $c_i \equiv 0$, then (4.5) is always satisfied, since there is at least one $i$ such that $Y_i = 1$. In fact, it is easily seen that (4.5) is a necessary and sufficient condition for the modified maximum likelihood estimating equation (4.4) to have a solution.

REMARK 6. The use of upper and lower bounds $M$ and $m$ for the $a_k$ is explained in the preceding section. The requirement that the $c_k$ be bounded above by $1 - \delta_0$ is natural as the guessing parameter would never exceed 0.5. However, as indicated by one reviewer, there should be other constraints in real applications (e.g., we can not allow the algorithm to only select items with $a = M$ and $c = 0$).

Theorem 2 can now be extended to cover the sequential design as just described for the three-parameter logistic model.

THEOREM 3. *For the sequential design defined in this section, the modified maximum likelihood estimating equation (4.4) has, with probability 1, a unique solution for all large n. The solution is strongly consistent (i.e. $\hat{\theta}_n \to \theta$ a.s.). Furthermore, provided that*

$$\frac{1}{v_n} \sum_{k=1}^{n} \frac{a_k^2}{8(1-c_k)^2}[1 - 20c_k - 8c_k^2 + (1+8c_k)^{3/2}] \underset{P}{\to} 1$$

*for some nonrandom sequence $v_n$,*

(4.6) $$\sqrt{v_n}(\hat{\theta}_n - \theta) \to_{\mathcal{L}} N(0,1).$$

*The normalizing constant $v_n$ in (4.6) may be replaced by the estimated Fisher information*

(4.7) $$I_n(\hat{\theta}_n) = \sum_{k=1}^{n} \frac{(1-c_k)a_k^2[e^{a_k(\hat{\theta}_n - b_k)}]^2}{[c_k + e^{a_k(\hat{\theta}_n - b_k)}][1 + e^{a_k(\hat{\theta}_n - b_k)}]^2}.$$

PROOF. As in the proof of Theorem 2, define $G(t) = e^t/(1+e^t)$, $\bar{G}(t) = 1 - G(t)$ and $\mathcal{F}_k = \sigma\{Y_j, \hat{\theta}_j, a_{j+1}, c_{j+1}; j \le k\}$. Applying the martingale local convergence theorem, we have that, analogous to (3.24),

(4.8) $$\sum_{k=1}^{\infty} \frac{w_k[Y_k - c_k - (1-c_k)G(a_k(\theta - b_k))]}{\sum_{j=1}^{k} w_j^2[c_j + (1-c_j)G(a_j(\theta - b_j))](1-c_j)\bar{G}(a_j(\theta - b_j))}$$

converges a.s.,

where $w_k = a_k(1 + \sqrt{1+8c_k})/(2c_k + 1 + \sqrt{1+8c_k})$. A slight modification of the proof leading to (3.10) can be constructed to show that

(4.9) $$P\left(\sum_{k=1}^{\infty} w_k^2[c_k + (1-c_k)G(a_k(\theta - b_k))](1-c_k)\bar{G}(a_k(\theta - b_k)) < \infty\right) = 0.$$

To provide a sketch to the proof of (4.9), let $A_1$ denote the event inside the probability sign in (4.9). Then, on $A_1$, $\lim|\hat{\theta}_n| = \infty$. We next prove, by contradiction, that $\limsup \hat{\theta}_n = \infty$ is impossible. Suppose that $\limsup \hat{\theta}_n = \infty$. Then, there is a subsequence $n_j$ such that $\hat{\theta}_k \le \hat{\theta}_{n_j}$, $k \le n_j$.

By the definition of $\hat{\theta}_n$, for $n \ge k_0$,

(4.10) $$\sum_{k=1}^{n} w_k[Y_k - c_k - (1-c_k)G(a_k(\hat{\theta}_n - b_k))] \le 0,$$

with the equality holding if and only if $\sum_{k=1}^{n} w_k Y_k > \sum_{k=1}^{n} w_k c_k$. From (4.10), we have

$$\sum_{k=1}^{n} w_k(1-c_k)[G(a_k(\hat{\theta}_n - b_k)) - G(a_k(\theta - b_k))]$$

(4.11)
$$\geq \sum_{k=1}^{n} w_k[Y_k - c_k - (1 - c_k)G(a_k(\theta - b_k))],$$

which converges to a finite limit on $A_1$. However, we can easily see that (3.12) and (3.13) still hold here. But they imply that the left-hand side of (4.11) can be arbitrarily small, which is a contradiction. Thus, $\limsup \hat{\theta}_n < \infty$ on $A_1$. Similarly, $\liminf \hat{\theta}_n > -\infty$ on $A_1$. Hence, $A_1$ must be a null set, and (4.9) holds.

Now, by the Kronecker lemma, we get from (4.8) and (4.9) that

(4.12)
$$\frac{\sum_{k=1}^{n} w_k[Y_k - c_k - (1 - c_k)G(a_k(\theta - b_k))]}{\sum_{k=1}^{n} w_k^2[c_k + (1 - c_k)G(a_k(\theta - b_k))](1 - c_k)\bar{G}(a_k(\theta - b_k))} \to 0$$
a.s.

Furthermore, by the definition of $\hat{\theta}_n$, for $n$ large enough such that $r_n \leq \theta$,

$$\sum_{k=1}^{n} w_k Y_k \leq \sum_{k=1}^{n} w_k[c_k + (1 - c_k)G(a_k(\hat{\theta}_n - b_k))],$$

which is $\leq \sum_{k=1}^{n} w_k[c_k + (1 - c_k)G(a_k(\theta - b_k))]$ if $\hat{\theta}_n = r_n$. Therefore, (4.12) implies that

(4.13)
$$\frac{\sum_{k=1}^{n} w_k(1 - c_k)[G(a_k(\hat{\theta}_n - b_k)) - G(a_k(\theta - b_k))]}{\sum_{k=1}^{n} w_k^2[c_k + (1 - c_k)G(a_k(\theta - b_k))](1 - c_k)\bar{G}(a_k(\theta - b_k))} \to 0$$
a.s.,

which is analogous to (3.15).

By examining the derivation following (3.15), we see that the same argument can be used to show that $\lim |\hat{\theta}_n| < \infty$ a.s. In particular, this implies that, for all large $n$, $\hat{\theta}_n$ is the solution to (4.4). It also implies, together with (4.13), that $\hat{\theta}_n \to \theta$ a.s.

Finally, we can apply the Taylor series expansion to (4.4) to obtain asymptotic normality. The argument is exactly the same as that in the proof of Theorem 2. $\square$

**5. Discussion.** CAT has become a popular mode of educational assessment in the United States. Examples of large-scale applications include the Graduate Record Examination (GRE), the Graduate Management Admission Test (GMAT), the National Council of State Boards of Nursing and the Armed Services Vocational Aptitude Battery (ASVAB). The most important component in a CAT is the item selection procedure that is used to select items during the course of the test. To date the most commonly used item selection procedure is the maximum Fisher information method. The motivation for maximizing the Fisher information is to make the trait estimator the most efficient. This can be achieved by recursively estimating $\theta$ with current available data and assigning further items adaptively. However, it is necessary to establish the corresponding theoretical properties for the maximum information approach.

The main objective of this paper is to tackle the sequential design and related convergence problems arising from the inherent mechanism of adaptive testing. It is clear that the logistic item response theory models are natural choices for CAT. We showed that, for the Rasch model, the usual plug-in adaptive design anchored in the current maximum likelihood estimator of the ability parameter converges to the optimal limit, and is therefore asymptotically efficient; moreover, the rate of the convergence can be characterized by the asymptotic normality of the maximum likelihood estimator. For the two-parameter logistic model, a similar asymptotic theory was developed based on an additional parameter modeling assumption that the discrimination power is restricted to a compact interval. Examples were given to illustrate that such restriction is necessary. As to the three-parameter logistic model, since the maximum likelihood estimating function is not generally a monotone function of the ability parameter, the maximum likelihood estimator may not be unique and, therefore, establishing convergence for the three-parameter logistic model is more complicated. Recognizing this potential problem, we proposed an asymptotically equivalent estimating function that is monotone in the ability parameter. Consistency and asymptotic normality were then proved for the adaptive design based on the modified maximum likelihood estimator.

The large scale implementation of CAT has created many interesting statistical issues in design, modeling and analysis. Our theory is established for the idealized setup that assumes existence of an infinite item pool. Even though, in reality, only finitely many items are available, the theory can still serve as a useful guidance to CAT practitioners as to how to choose an item selection strategy and how to design a simulation validation as well. In practice, simulation studies are always needed to help practitioners to evaluate the performance of their adaptive designs. According to the divergence examples created for the two-parameter logistic model, items with low discrimination should be used at the beginning of the test while items with high discrimination should be used at later stages. Therefore, a significant aspect of the new developments presented in this paper is to provide theoretical support to the item selection strategy of the a-stratified item selection method [Chang and Ying (1999)].

In order to design a good CAT algorithm, many complex controls are needed such as item exposure control and content balance. The item exposure rate for each item is defined as the ratio of the number of times the item is administered to the total number of examinees. Since CAT is designed to select the best items for each examinee, certain types of items tend to be most often selected by the computers, and many items are not selected at all, thereby making item exposure rates quite uneven. In addition, various non-statistical constraints need to be considered during item selection. Today's large-scale application of computer-based achievement tests and licensure exams has generated great challenges to test development. Maintaining content representation and other constraints is central to test defensibility and validity. Examples of the nonstatistical constraints include: a certain proportion of items should be selected from each content area, correct answers should fall approximately equally on options A, B, C and D, and a limited number of special items are allowed on a test, such as items with negative stems (e.g, "Which of the following choices is NOT true?"), just to name a few.

The $a$-stratified method was proposed with the objective of limiting the exposure of any given item by using that item at the most advantageous point in testing. The $a$-stratified method attempts to control item exposure by using less discriminating items early in the test, when estimation is least precise, and saving highly discriminating items until later stages, when finer gradations of estimation are required. One of the advantages of the $a$-stratified method is that it attempts to equalize the item exposure rates for all the items in the pool. Recently, methods of controlling content balance for the $a$-stratified method were proposed [see, e.g., van der Linden and Chang (2003), Yi and Chang (2003) and Cheng, Chang and Yi (2007)]. The advantages for using these methods are twofold: First, they allow the implementation of constraint on item selection in a-stratified adaptive testing; second, the constrained $a$-stratified methods may result in a set of theoretical advantages. It is evident that, by enforcing certain reasonable regularity conditions, the consistency results presented in this paper can be generalized to the constrained $a$-stratified methods, along with other reasonable item selection methods, such as the Bayesian item-selection criteria [see, e.g., van der van der Linden (1998)] and several Kullback–Leibler information based methods [see, e.g., Chang and Ying (1996)].

Similar procedures can be developed, and their properties can be established for other parametric item response theory models. A particularly useful class is the normal ogive models, in which the logistic link function is replaced by the normal distribution function. A minor technical complication in dealing with the normal ogive is that, even in the one-parameter case, the maximum likelihood estimating function is not monotone and may have multiple roots. But this complication may be dealt with by slightly modifying the estimating function, as we did for the three-parameter logistic model. The example presented in Remark 1 following Theorem 2 appears to be somewhat paradoxical in that the design is intended to increase efficiency by making the discrimination parameter large.

The example presented in Remark 1 following Theorem 2 also appears to be somewhat paradoxical, in that the design is intended to increase efficiency by making the discrimination parameter large. But a closer look at the design reveals that the inconsistency of $\hat{\theta}_n$ should be expected. This is because the amount of information at $\theta$, the true ability parameter, for the $k$th item may be extremely small when $b_k$ is not close enough to $\theta$ and $a_k$ is large. More specifically, when the magnitude of $a_k(\theta - b_k)$ is large, the Fisher information for the item is exponentially small, with the exponent proportional to $-|a_k(\theta - b_k)|$. Since under the normal circumstances, $b_k = \theta_{k-1}$ is about $O(k^{-1/2})$ away from $\theta$ [Chang and Stout (1993)], the choice $a_k = k^3$ effectively makes $|a_k(\theta - b_k)|$ very large. However, we still do not know if by choosing $a_k = o(\sqrt{k})$ it will be sufficient to guarantee the consistency of $\hat{\theta}_k$.

Finally, it should be pointed out that Mislevy and Wu (1996) and Mislevy and Chang (2000) showed that item selection in CAT leads to a design with missing data that are *missing at random* (MAR). Therefore, most of the standard theory for MLE holds from a missing data point of view.

## APPENDIX

PROOF OF INCONSISTENCY FOR EXAMPLE 1. On event $A$, we know that $\hat{\theta}_k, k \leq n_0$, are initialized so that $\hat{\theta}_1 = \hat{\theta}_0 - \varepsilon_0, \ldots, \hat{\theta}_{n_0} = \hat{\theta}_{n_0-1} - \varepsilon_0$ and $\hat{\theta}_k, k \geq n_0 + 1$, satisfy the maximum likelihood equations. We first claim that, on $A$,

$$(A.1) \qquad \hat{\theta}_{n_0+1} \leq \hat{\theta}_0 = \max_{0 \leq k \leq n_0} \hat{\theta}_k.$$

Recall that $b_k = \hat{\theta}_{k-1}$. So, (3.3) entails

$$(A.2) \qquad \frac{(n_0+1)^3}{1 + \exp[(n_0+1)^3(\hat{\theta}_{n_0+1} - \hat{\theta}_{n_0})]} = \sum_{k=1}^{n_0} \frac{k^3 \exp[k^3(\hat{\theta}_{n_0+1} - \hat{\theta}_{k-1})]}{1 + \exp[k^3(\hat{\theta}_{n_0+1} - \hat{\theta}_{k-1})]}.$$

Suppose that (A.1) does not hold. Then, $\hat{\theta}_{n_0+1} \geq \hat{\theta}_k, \ k \leq n_0$, implying

$$\text{left-hand side of (A.2)} \leq \tfrac{1}{2}(n_0+1)^3$$

and

$$\text{right-hand side of (A.2)} \geq \frac{1}{2} \sum_{k=1}^{n_0} k^3.$$

These two inequalities contradict (3.8). Thus, (A.1) holds.

Applying (3.3) to $\hat{\theta}_{n_0+2}$, we get

$$
\text{(A.3)} \quad \frac{(n_0+2)^3}{1+\exp[(n_0+2)^3(\hat{\theta}_{n_0+2}-\hat{\theta}_{n_0+1})]} + \frac{(n_0+1)^3}{1+\exp[(n_0+1)^3(\hat{\theta}_{n_0+2}-\hat{\theta}_{n_0})]}
$$

$$
= \sum_{k=1}^{n_0} \frac{k^3 \exp[k^3(\hat{\theta}_{n_0+2}-\hat{\theta}_{k-1})]}{1+\exp[k^3(\hat{\theta}_{n_0+2}-\hat{\theta}_{k-1})]}.
$$

Note that, on $A$, since $Y_k = 1, k \geq n_0+1$, $\hat{\theta}_{n_0+k}$ is increasing in $k$. From (A.3), we claim either $\hat{\theta}_{n_0+2} \leq \hat{\theta}_{n_0} + \varepsilon_0$ or

$$
\text{(A.4)} \qquad \hat{\theta}_{n_0+2} - \hat{\theta}_{n_0+1} \leq \frac{1}{(n_0+2)^2}.
$$

To prove this claim, suppose $\hat{\theta}_{n_0+2} \geq \hat{\theta}_{n_0} + \varepsilon_0 = \hat{\theta}_{n_0-1}$. Then,

$$
\text{(A.5)} \quad \frac{(n_0+1)^3}{1+\exp[(n_0+1)^3(\hat{\theta}_{n_0+2}-\hat{\theta}_{n_0})]} \leq \frac{(n_0+1)^3}{1+\exp[(n_0+1)^3 \varepsilon_0]} < \frac{1}{6},
$$

where the last inequality comes from (3.7). But $\hat{\theta}_{n_0+2} \geq \hat{\theta}_{n_0-1}$, implying that

$$
\text{(A.6)} \quad \sum_{k=1}^{n_0} \frac{k^3 \exp[k^3(\hat{\theta}_{n_0+2}-\hat{\theta}_{k-1})]}{1+\exp[k^3(\hat{\theta}_{n_0+2}-\hat{\theta}_{k-1})]} \geq \frac{n_0^3 \exp[n_0^3(\hat{\theta}_{n_0+2}-\hat{\theta}_{n_0-1})]}{1+\exp[n_0^3(\hat{\theta}_{n_0+2}-\hat{\theta}_{n_0-1})]} > \frac{1}{2}.
$$

Combining (A.5), (A.6) with (A.3), we have

$$
\frac{(n_0+2)^3}{1+\exp[(n_0+2)^3(\hat{\theta}_{n_0+2}-\hat{\theta}_{n_0+1})]} > \frac{1}{3},
$$

which, in conjunction with (3.6), entails (A.4).

Likewise, for $\hat{\theta}_{n_0+3}$, we claim one of the following must be true:

$$
\text{(A.7)} \qquad \hat{\theta}_{n_0+3} \leq \hat{\theta}_{n_0} + \varepsilon_0 = \hat{\theta}_{n_0-1},
$$

$$
\text{(A.8)} \qquad \hat{\theta}_{n_0+3} - \hat{\theta}_{n_0+1} \leq \frac{1}{(n_0+2)^2},
$$

$$
\text{(A.9)} \qquad \hat{\theta}_{n_0+3} - \hat{\theta}_{n_0+2} \leq \frac{1}{(n_0+3)^2}.
$$

To show this, suppose all of them are not true. Then the likelihood equation

$$
\text{(A.10)} \quad
\begin{aligned}
& \frac{(n_0+3)^3}{1+\exp[(n_0+3)^3(\hat{\theta}_{n_0+3}-\hat{\theta}_{n_0+2})]} \\
& + \frac{(n_0+2)^3}{1+\exp[(n_0+2)^3(\hat{\theta}_{n_0+3}-\hat{\theta}_{n_0+1})]} \\
& + \frac{(n_0+1)^3}{1+\exp[(n_0+1)^3(\hat{\theta}_{n_0+3}-\hat{\theta}_{n_0})]} = \sum_{k=1}^{n_0} \frac{k^3 \exp[k^3(\hat{\theta}_{n_0+3}-\hat{\theta}_{k-1})]}{1+\exp[k^3(\hat{\theta}_{n_0+3}-\hat{\theta}_{k-1})]}
\end{aligned}
$$

cannot hold since, in view of (3.6) and (3.7),

$$\text{left-hand side of (A.10)}$$

(A.11)
$$\leq \frac{(n_0+3)^3}{1+\exp(n_0+3)} + \frac{(n_0+2)^3}{1+\exp(n_0+2)} + \frac{1}{6} < \frac{1}{2}.$$

Furthermore,

(A.12)  $$\text{right-hand side of (A.10)} \geq \frac{n_0^3 \exp[n_0^3(\hat{\theta}_{n_0+3} - \hat{\theta}_{n_0-1})]}{1+\exp[n_0^3(\hat{\theta}_{n_0+3} - \hat{\theta}_{n_0-1})]} > \frac{1}{2}.$$

From (A.11) and (A.12), we obtain the desired contradiction and therefore the claim that one of (A.7)–(A.9) must hold is true.

In view of the preceding derivations, we have, for $k = 1, 2, 3$,

(A.13)
$$\hat{\theta}_{n_0+k} \leq \hat{\theta}_0 + \sum_{j=n_0+1}^{n_0+k} \frac{1}{j^2}.$$

We now apply the mathematical induction to show that (A.13) holds for every $k$. Suppose it is true for $k \leq j$. We claim that one of the following must hold:

(A.14)
$$\hat{\theta}_{n_0+j+1} \leq \hat{\theta}_{n_0} + \varepsilon_0,$$

(A.15)
$$\hat{\theta}_{n_0+j+1} - \hat{\theta}_{n_0+k} \leq \frac{1}{(n_0+k+1)^2} \qquad \text{for some } k \leq j.$$

This can be proved by showing that if none of the above inequalities holds, then the likelihood equation for $\hat{\theta}_{n_0+j+1}$ implies

$$\sum_{k=n_0+2}^{n_0+j+1} \frac{k^3}{1+\exp(k)} > \frac{1}{3},$$

which is a contradiction to (3.6). Clearly (A.14) or (A.15) and the induction assumption imply that (A.13) holds with $k = j+1$. Hence, (A.13) holds for every $k$ on event $A$. Thus, on $A$,

(A.16)
$$\begin{aligned}
\limsup_{n \to \infty} \hat{\theta}_n &< \hat{\theta}_0 + \sum_{j=n_0+1}^{n_0+k} \frac{1}{j^2} \\
&< \theta - 1 - \frac{\pi^2}{6} + \sum_{j=n_0+1}^{\infty} \frac{1}{j^2} \\
&< \theta - 1.
\end{aligned}$$
$$\square$$

**Acknowledgments.** The authors would like to thank two associate editors and three referees for their valuable comments.

DEPARTMENTS OF EDUCATIONAL PSYCHOLOGY
    AND PSYCHOLOGY
MC708
UNIVERSITY OF ILLINOIS
1310 S. SIXTH STREET
CHAMPAIGN, ILLINOIS 61820
USA
E-MAIL: hhchang@illinois.edu

DEPARTMENT OF STATISTICS
MC4690
COLUMBIA UNIVERSITY
1255 AMSTERDAM AVENUE
NEW YORK, NEW YORK 10027
USA
E-MAIL: zying@stat.columbia.edu